\def\cA{{\cal A}}
\def\cB{{\cal B}}
\def\cC{{\cal C}}
\def\cD{{\cal D}}

\def\cG{{\cal G}}

\def\cI{{\cal I}}

\def\cT{{\cal T}}

\def\pf{{\hfill$\Box$}}

\def\proof{{\noindent {\it Proof}.\ \ }}
\def\ra{{\rangle}}
\def\la{{\langle}}

\def\l({{\left(}}
\def\r){{\right)}}

\def\({{\Biggl(}}
\def\){{\Biggr)}}
\def\[{{\Biggl[}}
\def\]{{\Biggr]}}

\documentclass[12pt]{article}

\newtheorem{thm}{Theorem}
\newtheorem{clm}[thm]{Claim}
\newtheorem{cnj}[thm]{Conjecture}

\newtheorem{fct}[thm]{Fact}
\newtheorem{lem}[thm]{Lemma}

\usepackage{graphicx,amssymb,amsmath,verbatim}

\begin{document}

%
%
\title{On Universal Cycles for Multisets}

\author{
Glenn Hurlbert\thanks{
    Department of Mathematics and Statistics,
    Arizona State University,
    Tempe, Arizona 85287-1804,
    \texttt{hurlbert@asu.edu}
    }\\
Tobias Johnson\thanks{
    Department of Mathematics,
    Yale University,
    New Haven, CT 06520,
    \texttt{tobias.l.johnson@gmail.com}.
    Research supported in part by NSF grant 0552730.
    }\\
Joshua Zahl\thanks{
    Department of Mathematics,
    Caltech,
    Pasadena, CA 91125,
    \texttt{jzahl@zahl.ca}.
    Research supported in part by NSF grant 0552730.
    }\\
}
\date{February 18, 2008}
\maketitle


\vspace{0.5 in}

%
%
\begin{abstract}
A {\it Universal Cycle} for $t$-multisets of $[n]=\{1,\ldots,n\}$
is a cyclic sequence of $\binom{n+t-1}{t}$ integers from $[n]$
with the property that each $t$-multiset of $[n]$ appears exactly
once consecutively in the sequence.  For such a sequence to exist
it is necessary that $n$ divides $\binom{n+t-1}{t}$, and it is
reasonable to conjecture that this condition is sufficient for
large enough $n$ in terms of $t$.  We prove the conjecture
completely for $t\in\{2,3\}$ and partially for $t\in\{4,6\}$.
These results also support a positive answer to a question of Knuth.
\end{abstract}

\newpage

%
%
\section{Introduction}\label{Intro}

Problem 109 in Section 7.2.1.3 of Donald Knuth's \emph{The Art of
Computer Programming} \cite{Knu} lists the following sequence.
\begin{equation}\label{3mc5}
1112335\ 2223441\ 3334552\ 4445113\ 5551224
\end{equation}
The noteworthy property of this sequence is discovered by listing
every consecutive triple, including the two formed by wrapping the
sequence cyclically. Here we have
$$111, 112, 123,\ldots, 224, 241, 411\ .$$
Not only are the triples distinct, but they are still distinct
when considered as unordered multisets. In fact, each 3-multiset
of $[5]=\{1,2,3,4,5\}$ appears exactly once. Such a sequence is
called a {\it Universal Cycle for $3$-multisets of $[5]$}.  In
this paper we will use the shortened term {\it Mcycle}, and in
particular {\it$t$-Mcycle}, to refer to a universal cycle for
$t$-multisets of $[n]$, just as the term {\it Ucycle} has come to
refer to a Universal Cycle for $t$-subsets of $[n]$. Universal
Cycles for a wide range of combinatorial structures were
introduced in \cite{CDG}, and \cite{Hur} contains the most
up-to-date knowledge on Ucycles for subsets. In
this work we concern ourselves with $t$-Mcycles for
$t\in\{3,4,6\}$. Note that any permutation of $[n]$ is a 1-Mcycle,
and any eulerian circuit (which exists if and only if $n$ is odd)
of the complete graph (with loops) on $n$ vertices
is a 2-Mcycle. The main conjecture is the following.

\begin{cnj}\label{bigN}
For $t$ large enough in terms of $n$, Universal Cycles for
$t$-multisets of $[n]$ exist if and only if $n$ divides
$\binom{n+t-1}{t}$.
\end{cnj}

That the condition above is necessary follows from the fact that
each symbol is in the same number of multisets and hence must
appear equally often in the cycle. Our preceding comments indicate
that the conjecture is true for $t\in\{1,2\}$. In Section
\ref{GenPf} we prove the following theorem.

\begin{thm}\label{GenThm}
Let $n_0(3)=4,\ n_0(4)=5$ and $n_0(6)=11$. Then, for
$t\in\{3,4,6\}$ and $n\ge n_0(t)$, Mcycles for $t$-multisets of
$[n]$ exist whenever $n$ is relatively prime to $t$.
\end{thm}

This theorem verifies the conjecture for $t=3$, but leaves open
the case $n\equiv 2\mod 4$ for $t=4$ (the case $n\equiv 0\mod 4$
doesn't satisfy the necessary condition) and many cases for $t=6$.

More to the point, Knuth suggests in his solution to Problem 109 
that perhaps a $t$-Mcycle of $[n]$ exists if and only if a $t$-Ucycle
of $[n+t]$ exists (since then the two necessary conditions coincide).
Our proof of Theorem \ref{GenThm} sheds light on a correspondence
in this direction.

Because of the difficulty of extending these results for other
values of $t$, it is worth considering other methods of
construction.  
Sections \ref{IndPf} and \ref{ConvPf} re-prove the special case of 
$t=3$ from Theorem \ref{GenThm} using different methods that show
promise of being extended to larger $t$.
In Section \ref{IndPf} we describe an inductive technique that
is an extension of a method developed by Anant Godbole and
presented at the 2004 Banff conference on Generalizations of de~Bruijn 
Cycles and Gray Codes.  In Section \ref{ConvPf} we outline a
technique to convert $3$-Ucycles into $3$-Mcycles on the same ground
set.  While this technique is limited by the fact that the ground set 
remains unchanged, and thus divisibility considerations will make it 
difficult to generalize for $t>3$, it is a step in the right direction.

%
%
\section{Proof of Theorem \ref{GenThm}}\label{GenPf}

%
%
\subsection{Definitions}\label{Defs}

What is apparent in the sequence (\ref{3mc5}) is that each block has the same 
difference sequence modulo 5, namely $0011022$.
This is a key property of our constructions, as it has been in the construction
of all Ucycles to date.
This motivates the following definitions (the same terminology as in \cite{Hur}
except that ``form'' replaces ``$d$-set'').

Let $S=\{s_1,\ldots,s_t\}$, $s_i\le s_{i+1}$, be a $t$-multiset of
$[n]$. Define its {\it form}, $F(S)=(f_1, \ldots,f_t)$ by
$f_i=s_{i+1}-s_i$, where indices are computed modulo $t$ and
arithmetic is performed modulo $n$ (where $n$ is used in place of
$0$ as the modular representative). Two forms are {\it equivalent}
if one is a cyclic permutation of the other. Two forms belong to
the same {\it class} whenever one is any permutation of the other.
For example, with $t=5$ and $n=30$, each of the 30 multisets
$$\{1,11,21,21,1\}, \{2,12,22,22,2\},\ldots, \{30,10,20,20,30\}$$ 
belongs to the form $(10,10,0,10,0)$.  Also, the two forms
$$(10,10,10,0,0) {\rm\ and\ } (10,10,0,10,0)$$ 
make up the class $[10,10,10,0,0]$.
We will maintain the use of braces, parentheses, and brackets in order to
distinguish the various objects from one another.

%
%
\subsection{Basic Method}\label{Meth}

The main idea of the proof is to build a transition digraph whose edges
correspond to the set of all forms.
In the case that no bad patterns (to be defined in the next section) exist,
we show that the digraph is eulerian, thereby listing all forms by traversing
an eulerian circuit.
From there we repeat the circuit $n$ times in order to list all sets.
Note that this repetition is what gives rise to the block structure found
in the 3-Mcycle (\ref{3mc5}).

We begin by choosing a representative for each class. We
distinguish one of the coordinates of the form
$(f_1,\ldots,f_{t-1},f_t)$ (because of the equivalence among its
cyclic permutations we may assume it to be $f_t$) in the
representation $(f_1,\ldots,f_{t-1};f_t)$ so as to infer the
ordering $\{i,i+f_1,\ldots,i+f_1+\cdots+f_{t-1}\}$ of all its
multisets. This singled out coordinate is therefore unused in the
linear listing of each of these multisets.

Similarly, we may represent a class by $[f_1,\ldots,f_{t-1};f_t]$,
signifying that $f_t$ is distinguished (unused) in each of its
forms. It is important, then, that $f_t$ be unique in order to
avoid ambiguity---this is the reason for wanting good patterns.
For example, with $t=4$ and $n=7$, we can choose $[1,1,0;5] $ to
represent [1,1,0,5]. This determines the representations
$(1,1,0;5),(1,0,1;5)$ and $(0,1,1;5)$ of its three forms, of which
$(1,0,1;5)$ denotes the (ordered) forms $\{1,2,2,3\}$,
$\{2,3,3,4\}, \ldots$ and $\{7,1,1,2\}$.

Based on these choices we define the {\it transition graph}
$\cT_{n,t}$ as follows. We define the {\it prefix}, resp. {\it
suffix}, of the form representation $(f_1,\ldots,f_{t-1};f_t)$ to
be $((f_1,\ldots,f_{t-2}))$, resp. $((f_2, \ldots,f_{t-1}))$. Our
use of double parentheses denotes that these are the vertices in
the transition graph $\cT_{n,t}$ whose directed edges are
precisely the representations involved.

\begin{figure}
\centerline{\includegraphics[height=2.5in]{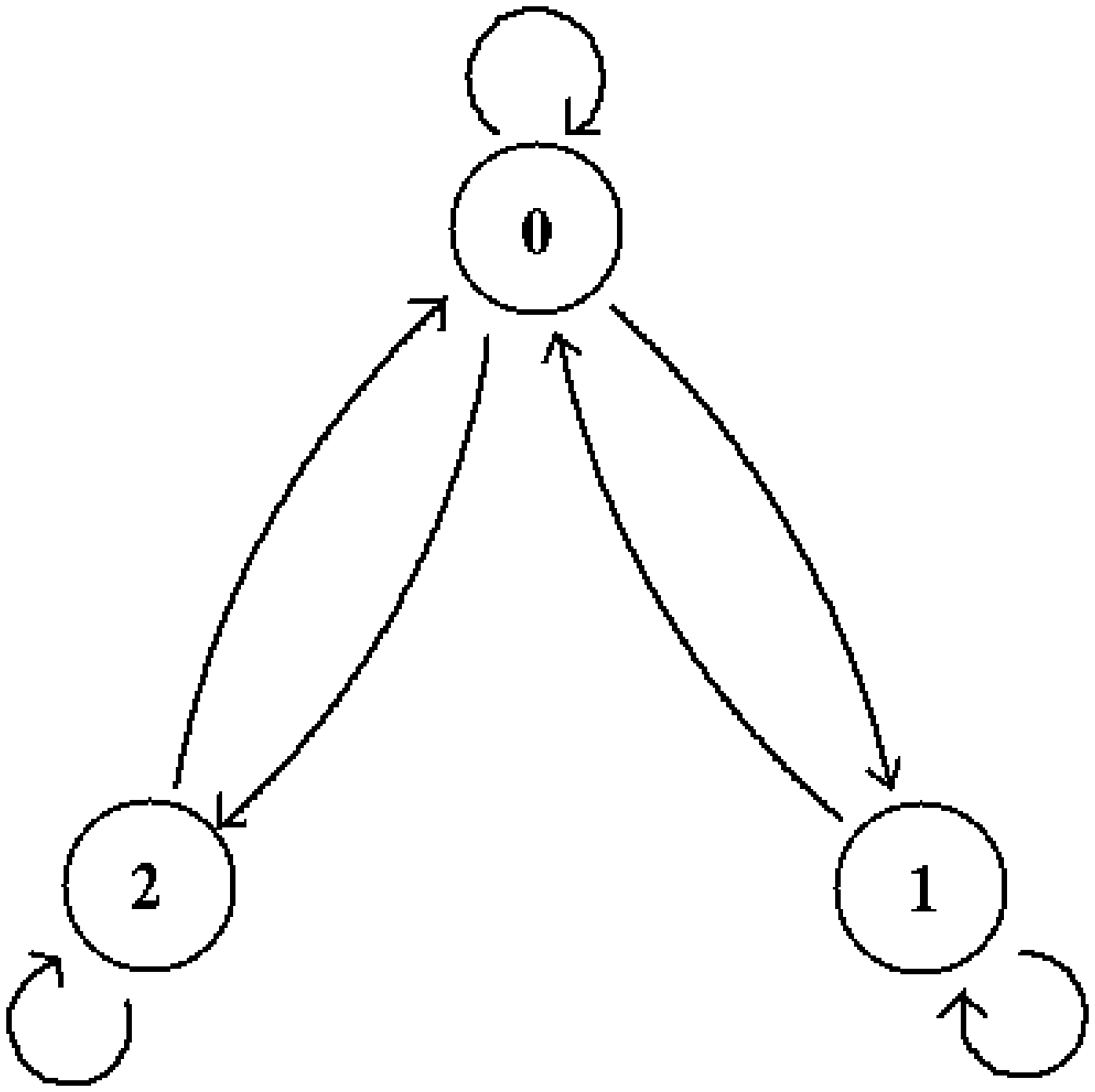}}
\caption{The transition graph $\cT_{5,3}$}\label{T53}
\end{figure}

For example, Figure \ref{T53} shows the transition graph
$\cT_{5,3}$, which was used to construct the 3-Mcycle (\ref{3mc5}).
The forms are
represented by $(0,0;5),(1,1;3)$, $(2,2;1),(0,1;4)$, $(1,0;4),
(0,2;3)$ and $(2,0;3)$. The form $(1,0;4)$ corresponds to the
directed edge $((1))\to((0))$, and so on. The eulerian circuit
0011022 corresponds to a listing of all forms and produces the
differences in the first block, 0001224, along with the first
digit, 1, of the next block. Since the sum $0+0+1+1+0+2+2\equiv
1\mod 5$, each block shifts by 1, and since 1 is relatively prime
to 5 each integer occurs as the starting digit of some block.
Hence, each 3-multiset of $[5]$ occurs exactly once. It turns out,
however (see the insertion technique in \cite{Hur}), that having 
the sum relatively prime to
$n$ is an unnecessary component in the general construction
process.

\begin{figure}
\centerline{\includegraphics[height=2.5in]{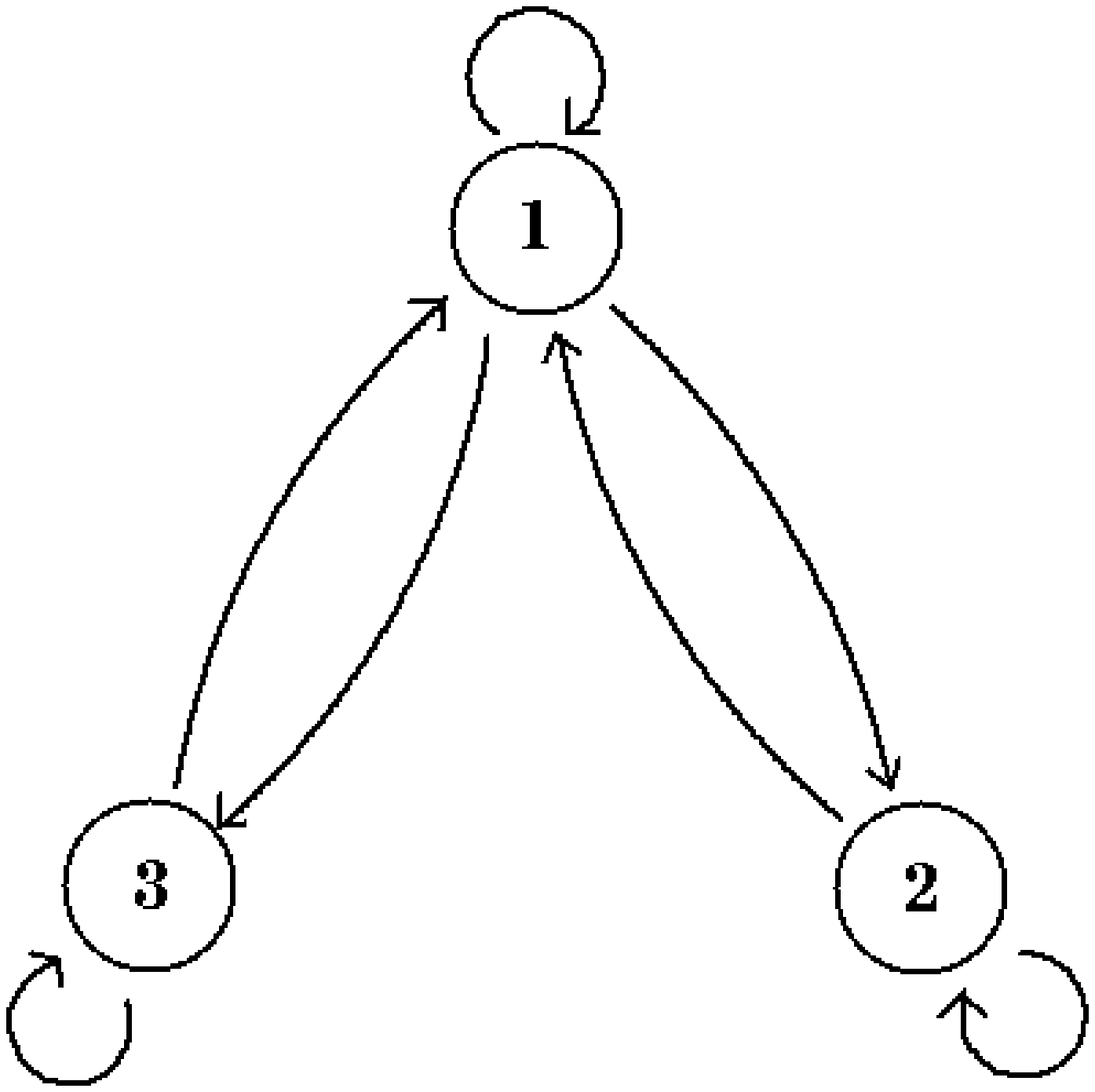}}
\caption{The transition graph $\cG_{8,3}$}\label{G83}
\end{figure}

In \cite{Hur} the analogous transition graph $\cG_{n,t}$ is
defined for $t$-subsets of $[n]$. For example, Figure \ref{G83}
shows the transition graph $\cG_{8,3}$, which can be used to
construct a Ucycle for 3-subsets of $[8]$. The forms in this case
are represented by $(1,1;6),(2,2;4)$, $(3,3;2),(1,2;5)$, $(2,1;5),
(1,3;4)$ and $(3,1;4)$. The eulerian circuit 1122133 generates the
Ucycle 
$$1235783\ 6782458\ 3457125\ 8124672\ 5671347\ 2346814\ 7813561\ 4568236\ ,$$
as above.

%
%
\subsection{Proof}\label{Pf}

In \cite{Hur} we find the following fact.

\begin{fct}\label{EulUcyc}
If $\cG_{n,t}$ is eulerian for some choice of representations of
classes, then there exists a Ucycle for $t$-subsets of $[n]$.
\end{fct}

The same arguments (as described in Section \ref{Meth}) that prove Fact
\ref{EulUcyc} yield the analogous result for Mcycles, which we therefore
state without proof.

\begin{lem}\label{EulMcyc}
If $\cT_{n,t}$ is eulerian for some choice of representations of
classes, then there exists an Mcycle for $t$-multisets of $[n]$.
\end{lem}

The key is that the obvious isomorphism between $\cT_{5,3}$ and
$\cG_{8,3}$ holds in general.

\begin{lem}\label{iso}
For every choice of representatives for the classes for
$t$-subsets of $[n+t]$ there exists choices of representatives for
the classes for $t$-multisets of $[n]$ so that the corresponding
transition graphs $\cG_{n+t,t}$ and $\cT_{n,t}$ are isomorphic.
\end{lem}
Of course, the theorem holds in reverse as well, with the roles of
$\cG$ and $\cT$ swapped, but we do not need that fact to prove
Theorem \ref{GenThm}. \medskip

\proof Simply shift every digit of every class and corresponding
representation down by one. Clearly $[f_1,\ldots,f_t]$ is a class
for $t$-subsets of $[n+t]$ if and only if $[f_1-1,\ldots,f_t-1]$
is a class for $t$-multisets of $[n]$.  The same can be said for
forms and for subsets/multisets.\pf

Now we borrow the final fact from \cite{Hur}.

\begin{fct}\label{EulG}
Let $n_0(3)=8$, $n_0(4)=9$, and $n_0(6)=17$. Then the transition
graph $\cG_{n,t}$ is eulerian for $t\in\{3,4,6\}$ and $n\ge
n_0(t)$ with $\gcd(n,t)=1$.
\end{fct}

In light of Lemma \ref{iso} and the knowledge that
$\gcd(n,t)=\gcd(n,n+t)$ we arrive at the following result.

\begin{lem}\label{EulT}
Let $n_0(3)=n_0(4)=5$ and $n_0(6)=11$. Then the transition graph
$\cT_{n,t}$ is eulerian for $t\in\{3,4,6\}$ and $n\ge n_0(t)$ with
$\gcd(n,t)=1$.
\end{lem}

The combination of Lemmas \ref{EulMcyc} and \ref{EulT} yields
Theorem \ref{GenPf}, with the exception that $n_0(3)=5$ instead of
4. However, a specific example for the case $t=3,\ n=4$ is given by
the sequence $S$ in (\ref{sequenceS})  at the end of Section \ref{IndPf},
below. This concludes the proof of Theorem \ref{GenThm}.
\pf
\medskip

Note that this result is weaker than Conjecture \ref{bigN} because
the relative primality condition replaces the divisibility
condition. 

To understand the limitations of this method we combine classes as follows.  
We see that the collection of all classes is the collection of all
unordered partitions of the integer $n$ into $t$ parts. Each class
defines a partition of $t$ according to the number of parts of the
same size. Two classes belong to the same {\it pattern} if they
define the same partition of $t$. We say that a pattern is {\it
good} if some part has size 1 and {\it bad} otherwise. 
By continuing the example from Section \ref{Defs} ($t=5, n=30$),
The 5 classes
$$[0,0,0,15,15], [2,2,2,12,12], [4,4,4,9,9], [8,8,8,3,3], {\rm\ and\ }
[10,10,10,0,0]$$
make up the pattern $\la 3,2\ra$.  
(Note that the class $[6,6,6,6,6]$ is skipped
from this sequence because it belongs to the pattern $\la 5\ra$.)
In this case there are 7 patterns in all:
$$\la 1,1,1,1,1\ra, \la 2,1,1,1\ra, \la 3,1,1\ra, \la 4,1\ra, \la 5\ra,
\la 2,2,1\ra, {\rm\ and\ } \la 3,2\ra\ ,$$
of which only $\la 5\ra$ and $\la 3,2\ra$ are bad. 

The underlying mathematics for Lemma \ref{EulT} comes
from the result of \cite{Hur} that no bad patterns exist if and
only if $t\in\{3,4,6\}$ and $\gcd(n,t)=1$.

%
%
\section{Inductive Construction for $t=3$}\label{IndPf}

As noted in the introduction, one of the main results of this work
is an inductive proof of Theorem \ref{GenThm} in the case that $t=3$.

\begin{thm}\label{IndThm}
For $n\ge 4$, $3$-Mcycles of $[n]$ exists provided $3\not|n$.
\end{thm}

\proof
For $t=3$, the condition $n\left|\binom{n+t-1}{t}\right.$ implies
that $n\equiv 1$ or $2$ (mod 3). We will show that for $n\geq 4$,
universal cycles on multisets exist whenever $n$ satisfies
$n\left|\binom{n+t-1}{t}\right.,$ i.e.
$n\left|\binom{n+2}{3}\right.$. We will prove this by induction on
$n$ as follows. We will start with a $3$-Mcycle on $[n-3]$ of the
form $stt^\prime uv$, where $stt^\prime uv$ is the concatenation
of the substrings $s,t,t^\prime, u,$ and $v$, where each of these
strings is a substring over the alphabet $[n-3]$ with specific
properties. From this string, we will construct a $3$-Mcycle on
$[n]$ of the form $STT^\prime UV$, where $S=st$, $T=t^\prime uv$,
$T^\prime$ is a cyclic permutation of $T$, and $U$ and $V$ are to
be described later.

Before describing the proof itself, we will define some
terminology that will be useful for describing universal cycles.
We say that a cyclic string $X=a_1a_2...a_k$ \emph{contains} the
multiset collection $\cI$ if $\cI=\big\{\{a_1,a_2,a_3\}$,
$\{a_2,a_3,a_4\}$, $\ldots$, $\{a_{k-2},a_{k-1},{a_k}\}$,
$\{a_{k-1},a_k,a_1\}$, $\{a_k,a_1,a_2\} \big\}$, where each of
these multisets must be distinct. Clearly $k=\binom{n+2}{3}$,
since this is the number of $3$-multisets on $[n]$.

For a string $X=a_1a_2...a_k$, we call the \emph{head} of $X$
the substring $a_1a_2$ and the \emph{tail} of $X$ the substring
$a_{k-1}a_k$.

Now, consider the collection of all 3-multisets over $[n]$. We
shall partition this collection into four subcollections. Let
$\cA$ be the collection of all 3-multisets over $[n-3]$, and let
$\cB$ be the collection of all $3$-multisets over $\{n-2,n-1,n\}$
and $[n-6]$ which contain at least one element from
$\{n-2,n-1,n\}$. Let $\cC$ be the collection of all 3-multisets
with one or two elements from $\{n-5,n-4,n-3\}$ and one or two
elements from $\{n-2,n-1,n\}$, and let $\cD$ be the collection of
all 3-multisets with one element from each of $[n-6]$,
$\{n-5,n-4,n-3\}$, and $\{n-2,n-1,n\}$. We can see that
$\cA,\cB,\cC,$ and $\cD$ are disjoint, and that their union is the
collection of all 3-multisets on $[n]$, as desired.

Now, let $S$ be a $3$-Mcycle on $[n-6]$, and since $1,1,1$ must
occur somewhere in $S$ and the beginning of $S$ is arbitrary, we
shall have $S$ begin with $1,1,1$. We shall also select $S$ so
that its tail is $n-6,n-7$. Thus $S$, when considered as a
cyclic string, contains all $3$-multisets over $[n-6]$, and when
considered as a non-cyclic string, contains all $3$-multisets
except $\{1,n-7,n-6\}$ and $\{1,1,n-7\}$. Let $T$ be a string over
$[n-3]$ such that $ST$---the concatenation of $S$ and $T$---is a
$3$--Mcycle over $[n-3]$. It is not clear that such a $T$ must
exist, but we shall find a specific example shortly. In the
example we will find, $T$ will begin with $1,1$ and will end with
$n-3,n-4$. Since $T$ begins with $1,1$, the string $ST$ contains
the multisets $\{1,n-7,n-6\},\ \{1,1,n-7\}$. We can see that the
cyclic string $ST$ contains all of the multisets in $\cA$, and
that when $ST$ is considered as a non-cyclic string, it contains
$\cA\backslash\big\{\{1,n-4,n-3\},\ \{1,1,n-4\}\big\}$. Now,
consider the string $T^\prime$ obtained by taking $T$ and
replacing each instance of $n-5$ by $n-2$, $n-4$ by $n-1$, and
$n-3$ by $n$. Since $T$ contained all multisets over $[n-3]$ which
contained at least one element from $\{n-5,n-4,n-3\}$, we have
that $T^\prime$ contains all multisets over $\{n-2,n-1,n\}$ and
$[n-6]$ which contain at least one element from $\{n-2,n-1,n\}$,
i.e. $T^\prime$ contains all the multisets in $\cB$. Since the
head of $T$ is $1,1$, the head of $T^\prime$ is also $1,1$,
and since $T$ ends with $n-3,n-4$, $T^\prime$ ends with $n,n-1$.
If we consider the cyclic string $STT^\prime$, we can see that
this string contains all the multisets in $\cA\cup \cB$, while the
non-cyclic version of this string is missing the multisets
$\{1,n-1,n\},\ \{1,1,n-1\}$.

For notational convenience, we will use the following assignments:
$a:=n-5,\ b:=n-4,\ c:=n-3,\ d:=n-2,\ e:=n-1,$ and $f:=n$. Now, we
will construct the strings $U$ and $V$. To do so, we shall
consider the case where $n$ is even and where $n$ is odd. For $n$
even, consider the following string:

\begin{eqnarray*}
V_e=&&\mathrm{be}(n-6)\mathrm{af}(n-7)\mathrm{be}(n-8)\mathrm{af}(n-9)...\mathrm{af}1\mathrm{be}\\
&&\ \ \ \mathrm{ad}(n-6)\mathrm{ce}(n-7)\mathrm{ad}(n-8)\mathrm{ce}(n-9)...\mathrm{ce}1\mathrm{ad}\\\
&&\ \ \ \ \ \
\mathrm{cf}(n-6)\mathrm{bd}(n-7)\mathrm{cf}(n-8)\mathrm{bd}(n-9)...\mathrm{bd}1\mathrm{cfe}.
\end{eqnarray*}
We can see that this string contains every multiset in $\cD$, as
well as the multisets $\{a,b,e\}$, $\{a,d,e\}$, $\{a,c,d\}$,
$\{c,d,f\},$ and $\{c,e,f\}$. Now, the following string (found
with the aid of a computer) contains all of the multisets in
$\cC\backslash \big\{\{a,b,e\}$, $\{a,d,e\}$, $\{a,c,d\}$,
$\{c,d,f\}$, $\{c,e,f\} \big\}$:
$$U_e=\mathrm{aaffc\phantom{1}aeebb\phantom{1}decec\phantom{1}bddcc\phantom{1}fbada\phantom{1}dfbf}.$$

Note that while the multisets $\{b,b,f\}$ and $\{b,e,f\}$ are not
present in the above string $U_e$, they are present in the
concatenation of $U_e$ with $V_e$. Similarly, while $U_e$ does not
contain $\{a,e,f\}$ and $\{a,a,e\}$, these multisets are present
in the concatenation of $T^\prime$ with $U_e$.

Now, we can see that the string $STT^\prime U_eV_e$ is a universal
cycle over $[n]$ because the non-cyclic string $STT^\prime$
contained all the multisets in
$\cA\cup\cB\backslash\big\{\{1,n-1,n\},\ \{1,1,n-1\} \big\}$, and
it is precisely the multisets $\{1,n-1,n\}$ and $\{1,1,n-1\}$
which are obtained by the wrap-around of the tail of $V_e$
with the head of $S$. The head and tail of the other
strings has been engineered so as to ensure that each multiset
occurs precisely once.

Now, consider the case where $n$ is odd. The corresponding strings
$V_o$ and $U_o$ are
\begin{eqnarray*}
V_o=&&\mathrm{be}(n-6)\mathrm{af}(n-7)\mathrm{be}(n-8)\mathrm{af}(n-9)...\mathrm{af}2\mathrm{be}\\
&&\ \ \ \mathrm{ad}(n-6)\mathrm{ce}(n-7)\mathrm{ad}(n-8)\mathrm{ce}(n-9)...\mathrm{ce}2\mathrm{ad}\\\
&&\ \ \ \ \ \
\mathrm{cf}(n-6)\mathrm{bd}(n-7)\mathrm{cf}(n-8)\mathrm{bd}(n-9)...\mathrm{bd}2\mathrm{cfe}.
\end{eqnarray*}
and
$$U_o=\mathrm{beb1f\phantom{1}abd1c\phantom{1}ffaae\phantom{1}cbfbf\phantom{1}dada1\phantom{1}eccfa\phantom{1}eecdc\phantom{1}dbd}.$$
The string $V_o$ contains the same multisets as $V_e$, with the
exception that $V_o$ does not contain the nine multisets
$\big\{\{1ad\},\{1ae\},...,\{1ce\},\{1cf\}\big\}$, and the string
$U_o$ contains the same multisets as $U_e$, with the exception
that it contains the additional nine multisets listed above. The
concatenation of $V_o$ and $U_o$ with the other strings works the
same way as their even counterparts.

This completes the induction proof, since the string $ST$ is a
$3$-Mcycle over $[n-3]$ (taking the place of $S$ in the previous
iteration of the induction), and the string $T^\prime UV$ extends
this $3$-Mcycle to $[n]$ (taking the place of $T$ in the previous
iteration of the induction). Also note that $T^\prime UV$ begins
with $1,1$ and ends with $n,n-1$, as required for the induction
hypothesis.

Thus, all that remains is the find a base case from which the
induction can proceed. A possible base case (there are many) for
$n=10$ is
\begin{eqnarray}\label{sequenceS}S&=&11144\phantom{1}42223\phantom{1}33121\phantom{1}24343\\
T&=&11522\phantom{1}63374\phantom{1}45166\phantom{1}27732\phantom{1}57366\phantom{1}77135\phantom{1}34641\phantom{1}71555\phantom{1}36127\phantom{1}42556\nonumber\\
&&66477\phantom{1}75526\phantom{1}4576\nonumber ,\end{eqnarray}
which would lead to
\begin{eqnarray*}
T^\prime&=&11822\phantom{1}93304\phantom{1}48199\phantom{1}20032\phantom{1}80399\phantom{1}00138\phantom{1}34941\phantom{1}01888\phantom{1}39120\phantom{1}42889\\
&&99400\phantom{1}08829\phantom{1}4809\\
U&=&55007\phantom{1}59966\phantom{1}89797\phantom{1}68877\phantom{1}06585\phantom{1}8060\\
V&=&69450\phantom{1}36925\phantom{1}01695\phantom{1}84793\phantom{1}58279\phantom{1}15870\phantom{1}46837\phantom{1}02681\phantom{1}709,
\end{eqnarray*}
Where ``0'' denotes 10 and the spacings have been added to
increase readability.

A possible base case for $n=11$ is
\begin{eqnarray*}S&=&11122\phantom{1}23114\phantom{1}22513\phantom{1}32444\phantom{1}33352\phantom{1}54541\phantom{1}43555\\ 
T&=&11657\phantom{1}43822\phantom{1}74468\phantom{1}54661\phantom{1}72736\phantom{1}18157\phantom{1}31888\phantom{1}77556\phantom{1}6688\phantom{1}57262\\
&&58536\phantom{1}21848\phantom{1}47776\phantom{1}41773\phantom{1}38826\phantom{1}
67836\phantom{1}36428\phantom{1}7,
\end{eqnarray*}
The corresponding strings $T^\prime,\ U,$ and $V$ can be found
using the method outlined above.
\pf

%
%
\section{Conversion Construction for $t=3$}\label{ConvPf}

\begin{thm}\label{ConvThm}
Any $3$-Ucycle of $[n]$ can be converted to a $3$-Mcycle of $[n]$.
\end{thm}

Before giving the proof, we introduce
two terms.  We call each element of $[n]$ a \emph{letter}, and
each $a_i$ in the Ucycle $X=a_1\ldots a_k$ a \emph{character}.  To
summarize, a 3-Mcycle on $[n]$ is made up of $\binom{n+t-1}{t}$
characters, each of which equals one of $n$ letters.

To demonstrate the proof's technique, we will first use an
argument similar to it to create 2-Mcycles from 2-Ucycles. We
start with this 2-Ucycle on $[5]$:
$$
1234513524.
$$
Then, we repeat the first instance of every letter to create the
following 2-Mcycle:
$$
112233445513524.
$$
The technique works because repeating a character $a_i$ as above
adds the multiset $\{a_i,a_i\}$ to the Ucycle and has no other
effect.

To use this technique on 3-Ucycles, we repeat not single
characters, but pairs of characters.  For example, changing
$$
\ldots a_{i-1}a_ia_{i+1}a_{i+2}\ldots
$$
to
$$
\ldots a_{i-1}a_ia_{i+1}a_ia_{i+1}a_{i+2} \ldots
$$
has only the effect of adding the 3-multisets
$\{a_i,a_i,a_{i+1}\}$ and $\{a_i,a_{i+1},a_{i+1}\}$ to the cycle.
In order to use this technique, we will need to know which
consecutive pairs of letters appear in a 3-Ucycle. For instance,
the following 3-Ucycle (generated using methods from \cite{Hur})
on $[8]$ contains every unordered pair of letters as consecutive
characters but $\{1,5\}$, $\{2,6\}$, $\{3,7\}$, and $\{4,8\}$:
$$1235783\ 6782458\ 3457125\ 8124672\ 5671347\ 2346814\ 7813561\
4568236\ .$$
This Ucycle is missing 4 pairs, which happens to be $n/2$.  This is no
coincidence: in fact, this is the most number of pairs that a
3-Ucycle can fail to contain.

\begin{clm}\label{missingpairs}
No two unordered pairs not appearing as consecutive characters in
a 3-Ucycle have a letter in common.  A 3-Ucycle can hence be
missing at most $n/2$ pairs of letters.
\end{clm}

\proof
Suppose that we have a 3-Ucycle that contains neither $a$ and $b$
as consecutive characters, nor $a$ and $c$ as consecutive
characters, where $a,b,c \in [n]$.  Then the 3-Ucycle does not
contain the 3-subset $abc$, because all permutations of $abc$ contain
either $a$ and $b$ consecutively, or $a$ and $c$ consecutively.
But this is a contradiction, as a 3-Ucycle by definition contains
all 3-subsets.

Hence, no two pairs of characters missing in the 3-Ucycle can have
a letter in common.  By the pigeonhole principle, the 3-Ucycle can
be missing at most $n/2$ pairs of letters.
$\diamondsuit$
\medskip

\noindent
{\it Proof of Theorem \ref{ConvThm}.}
Let $X$ be a 3-Ucycle on
$[n]$. Let $x_1,\ldots,x_n$ be a permutation of $[n]$ such that
\begin{itemize}
\item $x_1$ equals the first character in $X$. \item $x_n$ equals
the last character in $X$. \item If $x$ is even, the list
$\{x_1,x_2\}, \{x_3,x_4\}, \ldots, \{x_{n-1}x_n\}$ contains all
unordered pairs of letters not contained as consecutive characters
in $X$, which is possible by our lemma.  (If $X$ is missing
exactly $n/2$ pairs of letters, these pairs will be exactly the
pairs missing from $X$.  If $X$ is missing fewer than $n/2$ pairs
of letters, then the pairs consist of all missing pairs of
letters, plus the remaining letters paired arbitrarily.)

If $x$ is odd, one of the following lists contains all unordered
pairs of letters not contained as consecutive characters in $X$:
\begin{enumerate}
  \item $\{x_1,x_2\},\ldots,\{x_{n-2},x_{n-1}\}$
  \item $\{x_1,x_2\},\{x_4,x_5\},\{x_6,x_7\}\ldots,\{x_{n-1},x_{n}\}$
  \item $\{x_2,x_3\},\ldots,\{x_{n-1},x_n\}$
\end{enumerate}
This is possible by our lemma. There can be at most $(n-1)/2$
missing pairs of letters in $X$, and depending on whether $x_1$,
$x_n$, or both is a member of a missing pair, one of the above
lists can contain all the missing pairs. (As in the even case, it
does not present any problems if $X$ is missing fewer than
$(n-1)/2$ pairs.)
\end{itemize}
Construct $X^\prime$ by repeating the first instance of every unordered
pair of letters in $X$ except for $\{x_1,x_2\}, \{x_2,x_3\},$
$\ldots,$ $\{x_{n-1},x_n\},\{x_n,x_1\}$.  The cycle $X^\prime$ now
contains all multisets except \setlength\arraycolsep{1pt}
$$
\{x_1,x_1,x_1\},\ldots,\{x_n,x_n,x_n\}\\
$$
$$
\{x_1,x_1,x_2\},\{x_1,x_2,x_2\},\{x_2,x_2,x_3\},\{x_2,x_3,x_3\},\ldots,\{x_n,x_n,x_1\},\{x_n,x_1,x_1\}.
$$
Now, add the string $x_1x_1x_1x_2x_2x_2\ldots x_nx_nx_n$ to the
end of $X^\prime$ to create $X^{\prime\prime}$.  This provides
exactly the missing multisets, creating a 3-Mcycle.
\pf
\medskip

As a finall illustration, when $n=8$ we start with the 3-Ucycle
\begin{eqnarray*}
X & = & 1235783\ 6782458\ 3457125\ 8124672\\
&& 5671347\ 2346814\ 7813561\ 4568236\ .\\
\end{eqnarray*}
The 3-Ucycle $X$ does not contain the pairs $\{1,5\}$, $\{2,6\}$,
$\{3,7\}$, and $\{4,8\}$.  Hence, we set
\begin{eqnarray*}
x_1 & = & 1,\ x_2=5,\ x_3=3,\ x_4=7\\
x_5 & = & 4,\ x_6=8,\ x_7=2,\ x_8=6\ .
\end{eqnarray*}
Note that $x_1$ equals the first character of $X$, and $x_8$
equals the last.

Now, we repeat the first instance of every unordered pair except
for $\{1,5\}$, $\{5,3\}$, $\{3,7\}$, $\{7,4\}$, $\{4,8\}$,
$\{8,2\}$, $\{2,6\}$, and $\{6,1\}$.  (Note that four of these
pairs do not appear in $X$.  If some of these pairs actually did
appear in $X$, because $X$ was missing fewer than $n/2$ pairs of
letters, it would not affect the proof.)
\begin{eqnarray*}
X^\prime & = &12123235757878383\ 63676782424545858\ 3434571712525\
81812464672\\
&& 56567131347\ 2723468681414\ 7813561\ 4568236\ .
\end{eqnarray*}
Finally, we add the string $x_1x_1x_1\ldots x_nx_nx_n$ to complete
the Mcycle.
\begin{eqnarray*}
X^{\prime\prime} & = &12123235757878383\ 63676782424545858\
3434571712525\
81812464672\\
&& 56567131347\ 2723468681414\ 7813561\
4568236\\
&& 111555333777444888222666\ .
\end{eqnarray*}

%
%
\section{Remarks}\label{Remarks}

The proofs in Sections \ref{IndPf} and \ref{ConvPf} suggest
natural extensions to the $t=4$ and larger cases. Moreover, they
may prove useful by their introduction of new techniques for
approaching Ucycles. The section 3 proof is notable for its use of
induction, a technique which has not yet been used to create Ucycles.
This is especially promising in light of the many potential base
cases provided by Jackson \cite{Jac} for $t\le 11$. The section 4
proof, while it is tied to Ucycles, is not tied to any particular
approach for creating Ucycles, which is not true of the technique
in Section \ref{GenPf}.  Since the necessary condition for the
existence of Ucycles for $t$-subsets of $[n]$ is that $n$ divides
$\binom{n}{t}$, and since $\frac{1}{n}\binom{n+t-1}{t} =
\frac{1}{n+t}\binom{n+t}{t}$, we see that the condition for
$t$-Mcycles of $[n]$ is the same as the condition for $t$-Ucycles
of $[n+t]$. Thus it is reasonable to assume that some sort of
transformation between the two exists, as Knuth suggests.

For values of $n$ and $t$ for which Mcycles do exist, one
interesting question is how many Mcycles exist. Clearly each
Mcycles has $n!$ representations, since there are $n!$
permutations of $1,\ldots,n$. However, when searching for Mcycles
using a computer, vast numbers of \emph{distinct} (i.e. not
differing merely by a permutation of $1,\ldots,n$) Mcycles were
found. Currently, it is not clear whether $N(n,t)$, the number of
distinct Mcycles for a given value of $n$ and $t$, is a function
that can be approximated well.

%
%
\bibliographystyle{plain}
%

%
%
\end{document}